\documentclass{article}
\usepackage{amsmath,,amssymb,amsfonts}
\usepackage[american]{babel}
\usepackage{cite}
\usepackage{amsthm} 
\usepackage{indentfirst}
\usepackage[latin1]{inputenc}
\usepackage{geometry}
\usepackage{amsfonts}

\allowdisplaybreaks

 \newtheorem{theorem}{Theorem}[section] 
 \newtheorem{lemma}[theorem]{Lemma}

 \newtheorem{proposition}[theorem]{Proposition}
 \newtheorem{corollary}[theorem]{Corollary}

 \newenvironment{prf}{\begin{proof}[{\bf Proof}]}{\end{proof}}

\newcommand{\N}{{\mathbb{N}}}   
\newcommand{\Z}{{\mathbb{Z}}}   
\newcommand{\R}{{\mathbb{R}}}   
\newcommand{\Q}{{\mathbb{Q}}}   
\newcommand{\h}{{\mathcal{H}}}   
\newcommand{\hp}{{\mathbb{H}}}   
\newcommand{\C}{{\mathbb{C}}}   
\newcommand{\q}{{\mathbb{L}}}   
\newcommand{\U}{{\mathcal{U}}}   
\newcommand{\G}{{\mathcal{G}}}  
\newcommand{\n}{{\eta}}  

\newcommand{\para}{\par\vspace{.25cm}}

\begin{document}

 \title{Hyperbolic Unit Groups and Quaternion Algebras}

\author{S. O. Juriaans, I. B. S. Passi$^{*}$, A. C. Souza Filho}

\maketitle

\begin{abstract}
We classify  the quadratic  extensions $K=\Q [\sqrt{d}]$ and the
finite groups $G$ for which the  group ring $\mathfrak o_K[G]$ of
$G$ over the ring $\mathfrak o_K$ of integers of $K$ has the
property that the group $\mathcal U_1(\mathfrak o_K[G])$ of units of
augmentation 1 is hyperbolic. We also construct units in the
$\mathbb Z$-order $\h(\mathfrak o_K)$ of  the quaternion algebra
$\h(K)=\left(\frac{-1,\,-1}{K}\right)$, when  it is a division
algebra.   
  \end{abstract}
  

{\scriptsize{\textbf[keywords] {Hyperbolic Groups, Quaternion Algebras, Free Groups, Group Rings, Units }}}

\section{Introduction}

The finite groups $G$ for which the unit group $\mathcal U(\mathbb
Z[G])$ of the integral group ring $\mathbb Z[G]$ is hyperbolic, in
the sense of M. Gromov \cite{grmv}, have been characterized in
\cite{jpp}. The main aim of this paper is to examine the
hyperbolicity of the group $\mathcal U_1(\mathfrak o_K[G])$ of units
of augmentation 1 in the group ring $\mathfrak o_K[G]$ of $G$ over
the ring $\mathfrak o_K$ of integers of a quadratic extension
$K=\mathbb Q[\sqrt{d}]$ of the  field $\mathbb Q$ of rational numbers, where $d$
is a square-free integer $\neq 1$. Our main result (Theorem
\ref{mainresult}) provides a complete characterization of such group
rings  $\mathfrak o_K[G]$.

\para

In the integral case the hyperbolic unit groups are either finite,
hence have zero end, or have two or infinitely many ends (see
\cite[Theorem $I.8.32$]{mbah} and \cite{jpp}); in fact, in this
case, the hyperbolic boundary is either empty, or consists of two
points, or is a Cantor set. In particular, the hyperbolic boundary
is not a (connected) manifold. However, in the case we study here,
it turns out that when the unit group is hyperbolic and non-abelian
it has one end, and the hyperbolic boundary is a compact manifold of
constant positive curvature. (See Remark after Theorem  \ref{mainresult}.)

\para
Our investigation naturally leads us to study units in the order
$\h(\mathfrak o_K)$ of the standard quaternion algebra
$\h(K)=\left(\frac{-1,\,-1}{K}\right)$, when this algebra is a
division algebra. We  construct units, here called
Pell and  Gauss units, using  solutions of certain diophantine
quadratic equations. In particular, we exhibit units of norm $-1$
in $\h(\mathfrak o_{\mathbb Q[\sqrt{-7}]})$; this construction,
when combined with the deep work in \cite{cjlr}, helps  to provide a
set of generators for the full unit group $\mathcal U(\h(\mathfrak
o_{\mathbb Q[\sqrt{-7}]}))$.
\para
\hrule
\par\vspace{.25cm}
$^{*}$ {\small Senior Scientist, INSA.}

\newpage
\para
The work reported in this paper  corresponds to the first chapter of the third author's PhD thesis \cite{these}, where analogous questions about finite semi-groups, see \cite{ijsf},  and RA-loops, see \cite{jpmsf},  have also been studied.

\par\noindent
 
 
\section{Preliminaries}
Let $\varGamma$ be a finitely generated group with identity $e$ and
$S$ a finite symmetric set of generators of $\varGamma$, $e \not \in S$.
Consider the Cayley graph $\G = \G(\varGamma,S)$ of $\varGamma$ with
respect to the generating set $S$ and  $d=d_S$ the corresponding
metric (see \cite[chap. $1.1$]{mbah}). The induced metric on the
vertex set $\varGamma$ of $\G(\varGamma,S)$ is then the word metric:
for $\gamma_1,\,\gamma_2\in \varGamma$, $d(\gamma_1,\,\gamma_2)$
equals the least non-negative integer $n$ such that
$\gamma_1^{-1}\gamma_2=s_1s_2\ldots s_n,\ s_i\in S$ . Recall that in
a metric space $(X,\,d)$, the Gromov product $(y.z)_x$ of elements
$y,\,z\in X$ with respect to a given element $x\in X$ is defined to
be   $$(y \cdot z)_{x}=\frac{1}{2}(d(y,x)+d(z,x)-d(y,z)),
$$ and that the metric space $X$ is said to be {\it hyperbolic} if
there exists $\delta\ge 0$ such that
 for all $w,\,x,\,y,\,z \in X$, $$(x.y)_{w} \geq \min\{(x.z)_{w},\,(y.z)_{w}\}-\delta. $$ The group $\varGamma$ is said to be hyperbolic if the Cayley graph  $\G$ with the metric $d_S$ is a hyperbolic metric space.
This is a well-defined notion which depends only on the group
$\varGamma$, and is independent of the chosen generating set $S$
(see \cite{grmv}).
\para
 A map $f:X\to Y$ between topological spaces is said to be {\it proper} if 
 $f^{-1}(C)\subseteq X$ is compact whenever $C\subseteq Y$ is
 compact. For a
metric space $X$,  two proper maps (rays) $r_1,\,r_2:[0,\infty[
\longrightarrow X$  are defined to be equivalent if, for  each
compact set $C \subset X$, there exists $n \in \N$ such that $r_i([n, \infty[), i=1,\,2,$ are in the
same path component of $X \setminus C$. Denote by $end(r)$ the equivalence class of the ray $r$,  by $End(X)$  the set of the equivalence
classes $end(r)$, and by  $|End(X)|$ the cardinality of the set
$End(X)$. The cardinality  $|End(\G,\,d_S)|$ for the Cayley graph
$(\G,\,d_S)$ of $\varGamma$ does not depend on the generating set
$S$;  we thus have the notion of the number of ends of the finitely
generated group $\varGamma$  (see  \cite{mbah},\cite{grmv}).

\para
 We next recall some standard results from the theory of hyperbolic groups:
 \para
\begin{enumerate}

\item  \label{lgrm} Let $\varGamma$ be a group.
If $\varGamma$ is hyperbolic, then $\Z^{2} {\not \hookrightarrow}
\varGamma$, where $\Z^2$ denotes the free Abelian group of rank 2.\cite[Corollary III.$\Gamma\,3.10(2)$]{mbah}
\\
\item An infinite hyperbolic group contains an element of infinite order.\cite[Proposition III.$\Gamma\,2.22$]{mbah}\\
\item  \label{hppr}
If $\varGamma$ is hyperbolic, then there exists $n=n(\varGamma) \in
\N$ such that $|H|\leq n$ for every  torsion subgroup
$H<\varGamma$. \cite[Theorem III.$\Gamma\,3.2$]{mbah} and \cite[Chapter 8, Corollaire 36]{gh}
\end{enumerate}
These results will be used freely in the sequel. In view of (1)
above, the following  observation is quite useful.

\para
\begin{lemma} \label{hpli}
Let $A$ be a unital ring whose aditive group is torsion free, and let
$\theta_{1},\,\theta_{2} \in A$ be two  $2$-nilpotent commuting
elements  which are  $\Z$-linearly independent. Then $\U(A)$ contains a subgroup isomorphic to  $\Z^2$.
\end{lemma}

\begin{prf}
Set $u=1+\theta_{1}$ and $v=1+\theta_{2}$. It is clear that  $u,\,v
\in \U (A)$
 and both have infinite order. If $1 \neq w \in \langle u \rangle \cap
\langle v \rangle$ then there exists $i,\,j \in \Z \setminus
\{0\}$, such that, $u^{i}=w=v^{j}$. Since $u^{i}=1+i \theta_{1}$ and
$v^{j}=1+j \theta_{2}$, it follows that  $i \theta_{1}-j
\theta_{2}=0$ and hence $\{ \theta_{1},\, \theta_{2} \}$ is
$\Z$-linearly dependent, a contradiction. Hence $\mathbb Z^2\simeq
\langle u,\,v\rangle\subseteq\mathcal U(A).$
\end{prf}

\para
 Let  $C_{n}$ denote the cyclic group of order $n$, $S_3$ the
 symmetric group of degree 3, $D_4$ the
dihedral group of order $8$, and $Q_{12}$ the split extension
 $ C_3 \rtimes C_4$. Let $K$ be an algebraic number field and $\mathfrak o_K$
its ring of integers. The analysis of the implication for torsion
subgroups $G$ of a hyperbolic unit group $\mathcal U(\mathbb
Z[\varGamma])$ leading to  \cite[Theorem $3$]{jpp}  is easily seen to
remain valid for torsion subgroups of hyperbolic unit groups
$\mathcal U( \mathfrak o_K[\varGamma])$. We thus have the following:

\para

\begin{theorem}\label{tjpp}
A torsion group $G$  of a hyperbolic unit group $\mathcal
U(\mathfrak o_K[\varGamma])$ is isomorphic to one of the following
groups:
\para
\begin{enumerate}
\item $C_5,\, C_8,\, C_{12}$, an Abelian group of exponent dividing $4$
or $6$;\\
\item a Hamiltonian $2$-group;\\
\item $S_3,\, D_4,\, Q_{12},\, C_4 \rtimes C_4$.
\end{enumerate}
\end{theorem}

\para
We denote by  $\h(K)=(\frac{a,\,b}{K})$  the generalized quaternion
algebra over $K$:\linebreak
 $\h(K)=K[i,\,j:i^{2}=a,\,j^{2}=b,\,ji=-ij=:k].$ The set $\{1,\,i,\,j,\,k\}$ is
a $K$-basis of $\h(K)$. Such an algebra  is a totally definite
quaternion algebra if the field $K$ is totally real and $a,\,b$ are
totally negative. If $a,\,b \in \mathfrak o_K$, then the set
$\h(\mathfrak o_K)$, consisting  of the $\mathfrak o_K$-linear
combinations of  the elements  $1,\,i,\,j$ and $k$, is an $\mathfrak
o_K$-algebra. We denote by $N$ the norm map $\h(K)\to K$, sending 
$x=x_{1}+x_{i}i+x_{j}j+x_{k}k$
         to $N(x)=x_{1}^{2}-ax_{i}^{2}-bx_{j}^{2}+abx_{k}^{2}.$

\para                                       

Let $d\neq 1$ be a square-free integer, $K=\mathbb Q[\sqrt{d}]$.  Let us recall the basic facts
about the ring of integers $\mathfrak o_K$ (see, for example,
\cite{Janusz}, or \cite{LP}). Set $$\vartheta=\begin{cases}
\sqrt{d},\quad \text{if}\ d\equiv 2\ \text{or}\ 3\pmod 4\\
(1+\sqrt{d})/2,\quad  \text{if}\  d\equiv 1\pmod 4.
\end{cases}
$$
Then $\mathfrak o_K=\mathbb Z[\vartheta]$ and the elements $1,\
\vartheta$ constitute a $\mathbb Z$-basis of $\mathfrak o_K$. If
$d<0$, then
\begin{equation}\label{unit1}\mathcal U(\mathfrak o_K)=\begin{cases}
\{\pm 1,\ \pm \vartheta\},\quad \text{if}\ d=-1,\\
\{\pm 1,\ \pm \vartheta,\ \pm \vartheta^2\},\quad  \text{if} \ d=-3,\\
\{\pm 1\},\quad\text{otherwise}.
\end{cases}
\end{equation}
If $d>0$, then there exists a unique unit $\epsilon>1$, called the
{\it fundamental unit}, such that
\begin{equation}\label{unit2}\mathcal U(\mathfrak o_K)=\pm \langle
\epsilon\rangle.
\end{equation}
\para
We need the following:
\para
\begin{proposition} \label{salomao}
Let $K=\mathbb Q[\sqrt{d}]$, with $d\neq 1$ a square-free integer, be a
quadratic extension of $\mathbb Q$, and $u\in \mathcal U(\mathfrak
o_K).$ Then $u^i\equiv 1\pmod 2$, where \begin{equation*}
i=\begin{cases}1 \ \text{if}\  d\equiv 1\pmod 8,\\
2\ \text{if}\  d\equiv 2,\,3 \pmod 4,\\
3\ \text{if}\  d\equiv 5  \pmod 8.
\end{cases}
\end{equation*}
\end{proposition}
\para
\begin{prf}
The assertion follows immediately on considering the prime
factorization of the ideal $2\mathfrak o_K$, see \cite[Theorem $1$,
p.\,$236$]{BS}. \end{prf}


\par\vspace{.5cm}
\section{Abelian groups with hyperbolic unit groups} \label{seum}

\para
\begin{proposition}  \label{proli}
Let $R$ be a unitary commutative ring, $C_2=\langle g\rangle$. Then 
$u=a+(1-a)g$, $a\in R\setminus \{0,\,1\}$ is a non-trivial unit in
$\U_{1}(R[C_{2}])$
 if, and only if, $2a-1\in \mathcal U(R)$.
\end{proposition}

\para
\begin{prf}
 Let $C_{2}=\langle g \rangle$ and suppose that $u= a+(1-a)g,\, a\in R\setminus \{0,\,1\}$
 is a non-trivial unit in $R[C_2]$ having augmentation 1. Let $\rho:R[C_2]\to M_2(R)$ be the regular
 representation. Clearly $\rho(u)=\begin{pmatrix}
 a & 1-a \\ 1-a & a \end{pmatrix}.$ Since $u$ is a unit, it follows
that $2a-1= \det\rho(u)\in \U(R)$.
\para
 Conversely, let $ a \in  R
 \setminus \{0,1\}$ be such that $e= 2a-1 \in \U(R)$. It is then easy
 to see that $u=a+(1-a)g$ is a non-trivial unit in $R[C_2]$ with inverse
 $v=ae^{-1}+(1-ae^{-1})g$. \linebreak \cite[Proposition I]{lip} \end{prf}
\para
\begin{proposition}\label{trivial}
The unit group $\mathcal U(\mathfrak o_K[C_2])$ is trivial if, and
only if, $K=\mathbb Q $ or an imaginary quadratic extension of
$\mathbb Q$, i.e., $d<0$.
\end{proposition}\para
\begin{prf}
It is clear from the description  (\ref{unit1}) of the unit group of
$\mathfrak o_K$ that the equation

\begin{equation}\label{test}
2a-1=u,\ a\in \mathfrak o_k \setminus \{0,\,1\},\ u\in \mathfrak
\U(\mathfrak o_K)\end{equation} does not have a solution  when $K=\mathbb Q$ or
$d<0$.
\para
Suppose $d>1$ and $\epsilon$ is the fundamental unit in  $\mathfrak
o_K$. In this case we have  $\U(\mathfrak o_K)=\pm
\langle\epsilon\rangle$. By Proposition \ref{salomao}, $\epsilon
^{i}\in 1+2\mathfrak o_K$ for some $i\in \{1,\,2, \, 3\}$.
Consequently the equation  (\ref{test}) has a solution and so, by
Proposition \ref{proli},  $\mathcal U(\mathfrak o_K[C_2])$ is
non-trivial.
\end{prf}

\para
\begin{theorem} \label{tec}
Let $\mathfrak o_K$ be the ring of integers of a real quadratic
extension   $K=\mathbb Q[\sqrt{d}]$, $d>1$ a square-free integer,
$\epsilon>1$ the fundamental unit of $\mathfrak o_K$ and
$C_{2}=\langle g \rangle$.    Then  $$\U_{1}(\mathfrak o_K[C_{2}])
\cong \langle g \rangle \times \langle
\frac{1+\epsilon^{n}}{2}+\frac{1-\epsilon^{n}}{2}g \rangle \cong
C_{2} \times \Z,$$ where  $n$ is the  order of $\epsilon$ mod
$2\mathfrak o_K$.

\end{theorem}
\para
 \begin{prf}
 Let  $u\in \U_{1}(\mathfrak o_K[C_{2}])$ be a non-trivial unit. Then, there  exists $a \in
\mathfrak o_K$  such that, $2a-1=\pm\epsilon^{m}$ for some non-zero
integer $m$. Since $n$ is the order of   $\epsilon$ mod $2\mathfrak
o_K$, $m=nq$ with $q\in \mathbb Z$. We thus have

\par\vspace{.25cm}\noindent$ u= a+(1-a)g$
\par\vspace{.25cm}\noindent
 $\ =\frac{1\pm\epsilon^{m}}{2}+\frac{1\mp\epsilon^{m}}{2}g$
 \par\vspace{.25cm}\noindent
 $\ = \frac{1\pm\epsilon^{nq}}{2}+\frac{1\mp\epsilon^{nq}}{2}g$
 \par\vspace{.25cm}\noindent
 $\ = \left(\frac{1+\epsilon^{n}}{2}+\frac{1-\epsilon^{n}}{2}g\right)^q,$\
 or $ g\left(\frac{1+\epsilon^{n}}{2}+\frac{1-\epsilon^{n}}{2}g\right)^q.
 $

 \par\vspace{.25cm}\noindent
  Hence
  $\U_{1}(\mathfrak o_K[C_{2}]) \cong \langle g \rangle \times \langle \frac{1+\epsilon^{n}}{2}+\frac{1-\epsilon^{n}}{2}g \rangle \cong C_{2} \times \Z $.
\end{prf}
\para

As an immediate consequence of the preceding analysis, we  have:

\para
\begin{corollary}
If $K$ is a quadratic extension of $\ \mathbb Q$,  then $\mathcal
U_1(\mathfrak o_K[C_2])$ is a hyperbolic group.
\end{corollary}
\para
\begin{corollary}\label{klein} Let $G$ be a non-cyclic  elementary Abelian $2$-group.
Then $\U_{1}(\mathfrak o_K[G])$ is hyperbolic if, and only if,
$\mathfrak o_K$ is imaginary.
\end{corollary}

\para

\begin{prf} Suppose $\mathfrak o_K$ is real. Since $G$ is not cyclic, there exist $g, h \in G,\ g \neq \ h$, $o(g)=o(h)=2 $. By Theorem \ref{tec},
 $\U_{1}(\mathfrak o_K[\langle
  g\rangle])\cong C_{2} \times \Z \cong \U_{1}(\mathfrak o_K[\langle h\rangle])$. Since  $\langle
  g\rangle \cap \langle h\rangle=\{1\}$,  $\U_{1}(\mathfrak o_K[\langle g\rangle])\cap\, \U_{1}(\mathfrak o_K[\langle h\rangle])=\{1\}$.
  Therefore $\U_{1}(\mathfrak o_K)$ contains an Abelian group of rank $2$, so it is not hyperbolic.
 Conversely, if $\mathfrak o_K$ is imaginary, then, proceeding by induction on the order $|G|$ of $G$, we can conclude that  $\U_{1}(\mathfrak o_K[G])$
  is trivial, and   hence is hyperbolic.
\end{prf}



\para
   For an Abelian group $G$,  we denote by $r(G)$ its
  torsion-free rank. In order to study the hyperbolicity of $\U_1(\mathfrak o_K[G])$,
  it is enough to determine the torsion-free rank $r(\U_{1}(\mathfrak o_K[G]))$.
 Since $\U(\mathfrak o_K[G]) \cong \U(\mathfrak o_K) \times \U_{1}(\mathfrak o_K[G])$, we have
 $r(\U_{1}(\mathfrak o_K[G]))= r( \U(\mathfrak o_K[G]))- r( \U(\mathfrak o_K))$. If $K$ is an imaginary extension, then  $r(\U(\mathfrak o_K[G]))=r(\U_{1}(\mathfrak o_K[G]))$, whereas if $K$ is a real quadratic extension, then
   $r(\U( \mathfrak o_K))=1$, and therefore  $$r(\U_{1}(\mathfrak o_{K}[G]))= r(\U(\mathfrak o_K[ G]))-1.$$
 
We note that $$\mathbb Q[C_n]\cong \oplus \sum_{d|n}\mathbb
Q[\zeta_d],$$ where $\zeta_d$ is a primitive $d^{th}$ root of unity,
and therefore, for any algebraic number field $L$,
$$L[C_n]\cong \oplus \sum_{d|n}L\otimes_\Q \mathbb Q[\zeta_d].$$

We say that two groups are commensurable with each other when they contain finite index subgroups isomorphic to each other. Since
the unit group $\mathcal U(\mathfrak o_L[C_n])$ is commensurable
with $\mathcal U(\Lambda)$, where $\Lambda=\oplus
\sum_{d|n}\mathfrak o_{L\otimes \mathbb Q[\zeta_d]}$, we
essentially need to compute the torsion-free rank of $\mathfrak
o_{K\otimes \mathbb Q[\zeta_d]}$ for the needed cases.
\para
\begin{proposition} \label{rank} Let $K=\Q[\sqrt{d}]$, with $d$ a square-free integer $\neq 1$. The table below shows the   torsion-free rank of the groups $\U_{1}(\mathfrak o_K[C_{n}]),\  n \in\{2,\,3,\,4,\,5,\,6,\,8\}$.
$$\begin{tabular}{|c|l|c|l|}
\hline
$n$ & $r(\U_{1}(\mathfrak o_K[ C_{n}]))$&$n$ & $r(\U_{1}(\mathfrak o_K[ C_{n}]))$\\
\hline $2$ &  $\begin{array}{ll}
                                0&\textrm{ if } d<0\\
                                1& \textrm{ if } d>1
                                \end{array} $
& $3$ & $\begin{array}{ll}
                                1&\textrm{ if } d<0,d \neq -3 \\
                                0& \textrm{ if }d=-3\\
                1& \textrm{ if } d>1
                                \end{array} $ \\
\hline $4$ & $ \begin{array}{ll}
                         1&\textrm{ if } d<-1 \\
                         0& \textrm{ if } d=-1\\
                         2&\textrm{ if } d>1
                                \end{array}$
& $5$ & $\begin{array}{ll}
                      6&\textrm{ if } d<0 \\
                      2 & \textrm{ if } d=5\\
                      6&\textrm{ if } d\in \Z^{+}\setminus \{1,5\}
                                \end{array} $\\
\hline $6$ & $\begin{array}{ll}
                      2&\textrm{ if } d<-3 \\
                      0 & \textrm{ if } d=-3\\
                      3&\textrm{ if } d>1
                                \end{array} $
& $8$ & $\begin{array}{ll}
                             4&\textrm{ if } d<-1 \\
                             1& \textrm{ if } d=-1\\
                             4&\textrm{ if } d=2\\
                             5&\textrm{ if } d>2
                                \end{array} $\\
\hline
\end{tabular}$$
\end{proposition}
\par\vspace{.5cm}
In all  the cases, the computation is elementary and we omit the details.

\para
\begin{theorem} \label{trank}
If $K=\Q [\sqrt{d}]$,  with $ d$ a square-free integer $\neq 1$, then
\para
\begin{enumerate}
\item   $\U_{1}(\mathfrak o_K[C_{3}])$ is hyperbolic;\\
\item \label{hpc4}   $\U_{1}(\mathfrak o_K[C_{4}])$ is hyperbolic if, and only if,
  $d<0$;\\
\item for an Abelian group $G$ of exponent dividing $n>2$,
the group $\U_{1}(\mathfrak o_K[G])$ is hyperbolic if, and only if, $n=4$ and $d=-1$, or $n=6$ and $d=-3$ ;\\
\item   $\U_{1}(\mathfrak o_K [C_8])$ is hyperbolic if, and only if, $d=-1$;\\
\item   $\U_{1}(\mathfrak o_K[C_{5}])$ is not hyperbolic.
\end{enumerate}
\end{theorem}

\para
\begin{prf}$\empty$
The Proposition \ref{rank} gives us the torsion-free  rank
$$r:=r(\U_{1}(\mathfrak o_K[ C_{n}]))$$ for ${n \in
\{2,\,3,\,4,\,5,\,8\}}$. The group     $\U_{1}(\mathfrak o_K[
C_{n}])$ is hyperbolic if, and only if,  $r \in \{0,\,1\}$.  Thus, it only remains
to consider the case  $(3)$.
\para
Suppose $n=6$ and  $\mathcal U_1 (\mathfrak o_K[G])$ is hyperbolic.
We, hence, have $r \in \{0,1\}$. If  $G$ is cyclic, then, by   Proposition
\ref{rank}, we have $d=-3$. If $G$ is not cyclic,
 then $G \cong C_2^l \times C_3^m, \ l,\ m \geq 1$. Since  $\mathfrak o_K[C_3] \hookrightarrow \mathfrak o_K[G]$, it follows that $d=-3$.\para
 Conversely, if $n=6$ and $d=-3$ then, proceeding by induction on $|G|$, it can be proved that   $\U_{1}(\mathfrak o_K[G])$ is hyperbolic.
\para
 The case $n=4$ can be handled similarly. \end{prf}

\para
\begin{proposition}If $K=\Q [\sqrt{d}]$,  with $d$ square-free integer $\neq 1$,  then
  $\U_{1}(\mathfrak o_K[C_{12}])$ is not hyperbolic.
\end{proposition}

\para
\begin{prf} Since
 $\displaystyle K[ C_{12}] \cong K \otimes_{\Q}[\Q [C_{12}]] \cong K \otimes_{\Q}\Q [C_{3}\times C_{4}]) \cong K[C_{3}\times C_{4}]$,
  we have the immersions $\mathfrak o_K[C_{3}] \hookrightarrow \mathfrak o_K[C_{12}]$ and $ \mathfrak o_K[C_{4}] \hookrightarrow \mathfrak o_K[C_{12}]$.
  Therefore, $r(\U_{1}(\mathfrak o_K[C_{12}]))\geq r(\U_{1}(\mathfrak o_K[C_{3}]))+r(\U_{1}(\mathfrak o_K[C_{4}]))$.
 \para
   Suppose
    $\U_{1}(\mathfrak o_K[C_{12}])$ is hyperbolic. Then,  since $r(\U_{1}(\mathfrak o_K [C_{12}]))<2$, we have, by the Proposition \ref{rank}, $d \in \{-3,-1\}$.
 We also have $$K[C_{3}\times C_{4}]\cong (K[C_{3}])[C_{4}]\cong (K \oplus K[\sqrt{-3}])[C_{4}]
 \cong$$ $$ K[C_{4}] \oplus (K[\sqrt{-3}])[C_{4}] \cong 2K \oplus K[\sqrt{-1}] \oplus
 2K[\sqrt{-3}] \oplus K[\sqrt{-3}+\sqrt{-1}].$$
 \para
  Set
 $\q=\Q[\sqrt{-3}+\sqrt{-1}]$ and suppose
  $d=-3$. Then $\mathfrak o_K[C_{12}] \hookrightarrow 4\mathfrak o_K \oplus 2 \mathfrak o_{\q}$ and
 $r(\U(\mathfrak o_{\q}))=1$. Thus $r(\U(\mathfrak o_K[C_{12}]))=2$, and we have a contradiction.

 \par\vspace{.25cm} Analogously, for
 $d=-1$,  $\mathfrak o_K[C_{12}] \hookrightarrow 3\mathfrak o_K \oplus 3 \mathfrak o_{\q}$ and so
 $r(\U(\mathfrak o_K[C_{12}]))=3$. Since the extensions are non-real, we have that $r(\U_{1}(\mathfrak o_K
 [C_{12}]))=r(\U(\mathfrak o_K [C_{12}])) \geq 2$, and, hence, we again have a contradiction.
 \para
 We conclude that $\mathcal U_1(\mathfrak o_K[C_{12}])$ is not hyperbolic. \end{prf}


\par\vspace{.5cm}
\section{ Non-Abelian groups with hyperbolic unit groups}
\para
 Theorem \ref{tjpp} classifies the finite non-Abelian groups $G$ for which the unit group
$\U_{1}(\Z [G])$ is hyperbolic. These groups are:
 $S_{3},\,D_{4},\,Q_{12}, C_{4}\rtimes C_{4}$,  and the Hamiltonian  $2$-group, where \linebreak $Q_{12}=C_3 \rtimes C_4$, with $C_4$ acting non-trivially on $C_3$, and also on $C_4$
(see \cite{jpp}).

\para

  E. Jespers, in \cite{jsp}, classified the finite groups $G$ which have a normal
 non-Abelian free complement in $\U (\Z [G])$.
The group algebra $\Q [G]$ of these groups has at most one matrix
Wedderburn component which must be isomorphic to $M_{2}(\Q)$.
\para
\begin{lemma} \label{mhp}
Let $G$ be a group and $K$ a quadratic extension. If $M_{2}(K)$ is a
Wedderburn component of $K[G]$ then
 $\Z^{2} \hookrightarrow \U_{1} (\mathfrak o_K[G])$. In particular, $\U_{1} (\mathfrak o_K[G])$ is not   hyperbolic.
\end{lemma}

\begin{prf} The ring
$\varGamma=M_{2}(\mathfrak o_K)$ is a $\Z$-order in $M_{2}(K)$ and
  $$X=\{e_{12},\,e_{12}\sqrt{d}\} \subset \varGamma$$ is a set of commuting nilpotent elements of index $2$, where $e_{ij}$ denotes the elementary matrix. The set $\{1,\sqrt{d}\}$ is a linearly independent set over $\mathbb Q$, and hence so is $X$. Therefore, by  Lemma \ref{hpli},
  $\Z^2 \hookrightarrow \U_{1} (\varGamma) \subset \U_{1} (\mathfrak o_K[G])$, and so,  $\U_{1} (\mathfrak o_K[G])$ is not   hyperbolic.
\end{prf}

\para
\begin{corollary}
 If $\ G \in \{S_{3},\,D_{4},\,Q_{12},\,C_{4}\rtimes C_{4}\}$ then $\U_{1}(\mathfrak o_K[G])$
 is not hyperbolic.
\end{corollary}

\para
\begin{prf}
We have that $\displaystyle K [G] \cong K \otimes_{\Q}(\Q [G])$. For
each of the groups under consideration,    $M_{2} (\Q)$  is a
Wedderburn component of  $\Q [G]$; it therefore follows that
$M_{2}(K)$ is a Wedderburn component of
 $K [G]$. The preceding  lemma implies that   $\U_{1}(\mathfrak o_K[G])$ is not hyperbolic.
\end{prf}

\para
If $H$ is a non-Abelian Hamiltonian $2$-group,  then $H=E \times
Q_8$, where $E$ is an elementary Abelian $2$-group and $Q_8$ is the
quaternion group of order 8. Since $Q_{8}$ contains a cyclic
subgroup of order $4$, it follows, by  Theorem \ref{trank}, that {\it
if $\U_{1}(\mathfrak o_K [Q_{8}])$ is hyperbolic, then $\mathfrak
o_K$ is not real}.
\para
\begin{proposition}
If  $G$ is a Hamiltonian $2$-group of order greater than $8$, then
$\U_{1}(\mathfrak o_K[ G])$ is not hyperbolic.
\end{proposition}

\para
\begin{prf}
Let $G=E \times Q_{8}$ with  $E$ elementary Abelian of order
$2^{n}>1$. We then have  $K[G]=K[E \times
 Q_{8}] \cong
 \displaystyle K \otimes_{\Q}(\Q[E \times Q_{8}])\cong \displaystyle K \otimes_{\Q}(\Q [E])[Q_{8}] \cong \displaystyle K \otimes_{\Q}(2^{n} \Q)[Q_{8}]\cong
 (2^{n}K)[Q_{8}]$. If $d=-1$ it is well known that $K Q_8$  has a Wedderburn component isomorphic to $M_2(K)$ and hence, by Lemma \ref{mhp}, $\U_{1}(\mathfrak o_K Q_8)$ is not hyperbolic. If $d<-1$, then, by Proposition \ref{rank},
 $r(\U_1(\mathfrak o_K[ C_{4}]))=1$. Since  $C_{4}$ is a subgroup of  $Q_{8}$, it follows that
 $\U_1((2^n\mathfrak o_K)[C_4])$ embeds into $\U_1(\mathfrak o_K[G])$. Thus, since  $\U_1(\prod_{2^n} \mathfrak o_K[C_4])$ has rank $2^n \geq 2$,  $\U_1(\mathfrak o_K[G])$ is not hyperbolic.
\end{prf}
\para In view of the above Proposition, it follows that $Q_8$ is the only Hamiltonian $2$-group   for which
  $\U_{1}(\mathfrak o_K[G])$ can  possibly be hyperbolic,  and in this case
 $\mathfrak o_K$ is the ring of integers  of an imaginary extension. By Lemma \ref{mhp},  $K[Q_8]$ can not have a matrix ring as a Wedderburn component.  Since $\Q [Q_{8}] \cong 4 \Q \oplus \h(\Q)$, we have
 $K [Q_{8}] \cong \displaystyle K \otimes_{\Q}(4\Q \oplus \h(\Q)) \cong 4K \oplus \h(K) $; hence  $K[Q_8]$ must be a direct sum of division rings, or equivalently, has no non-zero nilpotent elements. In particular, $\h (K)$ is a division ring.

\para

\begin{theorem}
Let $K=\mathbb Q[\sqrt{d}]$, with $d$ square-free integer $\neq 1$. Then $K[Q_8]$ is a direct sum of division rings if, and only if, one of the following holds:\\
\newcounter{nil}
\begin{list}{$(\roman{nil})$}{\usecounter{nil}}
\item $d\equiv 1\pmod 8$;\\
\item $d\equiv 2,\ or\ 3\pmod 4,\ or \ d\equiv 5 \pmod 8,\  and\  d>0$.
\end{list}
\end{theorem}
\para
\begin{prf}
The assertion follows from \cite[Theorem 2.3]{Arora} ;
\cite[Theorem 1, p.\,236]{BS}
 and \cite[Theorem $3.2$]{rajw}.
\end{prf}

\para
\begin{corollary}
If  $K=\Q[\sqrt{d}]$, where $d$ is a negative square-free integer,
then the group $\U_{1}(\mathfrak o_K[Q_{8}])$ is not hyperbolic if
${d {\not \equiv} 1 \pmod 8}$.
\end{corollary}


\para
Let $\hp:\C \times ]0, \infty[$ be the upper half-space model of
three-dimensional hyperbolic space and $Iso(\hp)$ its group of
isometries.   In the quaternion algebra $\h:=\h(-1,-1)$ over $\R$,
with its usual basis, we may identify   $\hp$
 with the subset  $  \{z+rj: z \in \C, r \in \R^{+}\}.$
The group  $PSL(2,\,\C)$  acts on  $\hp$ in the following way: 
$$\begin{array}{llll}\varphi:&PSL(2,\,\C) \times \hp &\longrightarrow &\hp\\
                             &(M,\,P) &\mapsto &\left ( \begin{array}{ll}a&b\\c&d \end{array}\right )P:{=}(aP+b)(cP+d)^{-1}, \end{array}$$
where $(cP+d)^{-1}$ is calculated in $\h$. Explicitly,
$MP=M(z+rj)=z^{*}+r^{*}j$,
 with $$z^{*}=\frac{(az+b)(\overline{c}\overline{z}+\overline{d})+a \overline{c}r^{2}}{|cz+d|^{2}+|c|^{2}r^{2}},\ \text{and}\  r^{*}=\frac{r}{|cz+d|^{2}+|c|^{2}r^{2}}.$$

\para

 Let $K$ be an algebraic  number field and $\mathfrak o_K$ its
ring of integers.
 Let $$SL_{1}(\h(\mathfrak o_K)):=\{x \in \h(\mathfrak o_K):N(x)=1\},$$ where
$N$ is the norm in $\h(K)$. Clearly  the groups   $\U(\h(\mathfrak
o_K))$ and $\U(\mathfrak o_K) \times SL_{1}(\h(\mathfrak
o_K))$ are commensurable. Consider the subfield  $F= K[i]
\subset \h(K)$ which is a maximal subfield in   $\h(K)$. The inner
automorphism
$\sigma$, $$\begin{array}{llll} \sigma:& \h(K) &\longrightarrow &\h(K)\\
                                             & x &\mapsto &jx j^{-1},
                                             \end{array}$$
fixes $F$. The algebra  $\h(K) =F \oplus Fj$ is a crossed product
 and    embeds into  $M_{2}(\C)$ as follows:
\begin{equation}\label{reph}\begin{array}{llll} \Psi:& \h(K) &\hookrightarrow &M_{2}(\C)\\
                                             & x+yj &\mapsto &\left (
                                             \begin{array}{ll}x&y\\-\sigma(y)&\sigma(x)
                                             \end{array}\right ).
                                             \end{array}  \qquad\qquad\qquad \end{equation}

\para

This embedding enables us to  view  $SL_{1}(\h(\mathfrak o_K))$ and
$SL_{1}(\h(K))$ as subgroups of $SL(2,\,\C)$ and hence
$SL_{1}(\h(K))$ acts on $\hp$.
\para
\begin{proposition}\label{hypq8}
Let $K=\mathbb Q[\sqrt{d}]$, $d\equiv 1\pmod 8$ a square-free
negative integer, and $\mathfrak o_K$ its ring of integers. Then
$\mathcal U(\h(\mathfrak o_K))$ and   $\mathcal U(\mathfrak
o_K[Q_8])$ are hyperbolic groups.
\end{proposition}
\para
\begin{prf} Observe that  $SL_{1}(\h(\mathfrak o_K))$ acts on the space $\hp$ and, hence, is a discrete subgroup of $SL_{2}(\C)$ (see \cite[Theorem $10.1.2$, p.\,$446$]{egm}). 
The quotient space 
 $Y:=\hp /SL_{1}(\h(\mathfrak o_K))$ is a Riemannian manifold of constant curvature $-1$ and,  since $\hp$ is simply connected, we have that $SL_{1}(\h(\mathfrak o_K))\cong\pi_{1}(Y)$. 
 Since  $d \equiv 1 \pmod{8}$, $\h(K)$ is
a division ring and, therefore, co-compact and $Y$ is compact  (see  \cite[Theorem $10.1.2$, item $(3)$]{egm}).  Hence $SL_{1}(\h(\mathfrak o_K))$
 is hyperbolic  (see \cite[ Example $2.25.5$]{bk}).
 Since $\U(\h(\mathfrak o_K))$  and  $\U(\mathfrak o_K) \times SL_{1}(\h(\mathfrak o_K))$
 are
 commensurable and $\U(\mathfrak o_K)=\{-1,\,1\}$, it follows that $\U(\h(\mathfrak o_K))$ is   hyperbolic. Since $\U(\mathfrak o_K[ Q_8])\cong \U(\mathfrak o_K) \times  \U(\mathfrak o_K) \times  \U(\mathfrak o_K) \times  \U(\mathfrak o_K) \times \U(\h(\mathfrak o_K))$
and $\U(\mathfrak o_K) \cong C_2$,  we conclude that
  $\U(\mathfrak o_K[ Q_8])$ is  hyperbolic.
\end{prf}
\para
Combining the results in the present and the preceding section,
we have the following  main result.
\para
\begin{theorem}\label{mainresult} Let $K=\Q[\sqrt{d}]$, with $d$ square-free integer $\neq 1$, and $G$ a finite group.
   Then $\U_{1}(\mathfrak o_K[ G])$ is
 hyperbolic if, and only if, $G$ is one of the groups listed  below
 and $\mathfrak o_K$ (or $K$) is determined by the corresponding value of $d$:\\

\begin{enumerate}
 \item $G \in \{ C_{2}, \,C_{3} \}$ and $d$ arbitrary;\\
 \item $G$ is an Abelian group of exponent dividing
  $n$ for:\\$n=2$ and $d<0$, or $n=4$ and $d=-1$, or $n=6$ and $d=-3$.\\
 \item $G=C_{4}$ and $d<0$.\\
  \item $G=C_8$ and $d=-1$.\\
  \item $G=Q_{8}$ and
   $d<0$ and ${d \equiv 1 \pmod{8}}$.
\end{enumerate}
\end{theorem}

\para

 \noindent
{\bf Remark.} If the group $\U(\mathfrak o_K[ Q_8])$ is hyperbolic then
the hyperbolic boundary   ${\partial(\U(\mathfrak o_K[ Q_8]))\cong
\mathbb S^2}$, the Euclidean sphere of dimension $2$,  and
$Ends(\U(\mathfrak o_K[ Q_8]))$ has one element (see \cite[Example $2.25.5$]{bk}). Note that if $\U (\Z [G])$ is an infinite
non-Abelian hyperbolic group, then $\partial(\U(\Z[ G]))$ is totally
disconnected and is a Cantor set. So, in this case, $\U(\Z [G])$ has
infinitely many ends and also is a virtually free group,
(\cite[Theorem $2$]{jpp} and \cite[\S $3$]{grmv}). However, if
$\U(\mathfrak o_K[G])$ is a non-Abelian hyperbolic group, then
$\U(\mathfrak o_K[G])$ is an infinite group which
 is not virtually free, has one end and $\partial(\U(\Z[ G]))$
is a smooth manifold.
\par\vspace{.5cm}

\section{Pell and  Gauss Units}
\para
 When the algebra $\h(K)$ is isomorphic to $M_{2}(K)$ it
is known how to construct the unit group of a $\Z$-order up to a
finite index. Nevertheless, if $\h(K)$ is a division ring, this is a
highly non-trivial task; see  \cite{cjlr}, for example. In this
section we study  construction of  units of $\U(\h(\mathfrak o_K))$
in the case when the quaternion algebra $\h(K)$ is a  division ring.

\para

In the sequel,  $K=\Q[\sqrt{-d}]$  is  an imaginary quadratic
extension with  $d$ a square-free integer congruent to $ 7\pmod 8$,
and $\mathfrak o_K$  the ring of  integers of the field $K$. Note
that $s(K)$, the {\it stufe} of $K$, is 4, the quaternion algebra
$\h(K)$ is a division ring and $\mathcal U(\mathfrak o_K)=\{\pm
1\}$. Thus, if $u=u_{1}+u_{i}i+u_{j}j+u_{k}k \in \U(\h(\mathfrak
o_K))$,  then its norm
 ${N(u)=u_{1}^{2}+u_{i}^{2}+u_{j}^{2}+u_{k}^{2}=\pm 1}$; furthermore, if any of the coefficients $u_{1},\,u_{i},\,u_{j},\,u_{k}$ is zero then $N(u)=1$, $s(K)$ being 4.
\para

The representation of $u$, given by (\ref{reph}), is
$$[u]:=\Psi(u)=\begin{pmatrix} u_1+u_ii &u_{j}+u_ki\\ -u_{j}+u_k i& u_1-u_ii\end{pmatrix}\in M_2(\mathbb C).$$
 \para
 Denote by  $\chi_u$  the characteristic polynomial
of $[u]$, and by $m_u$ its minimal polynomial. The degree
$\partial(\chi_u)$ of $\chi_u$ is $2$ and therefore $\partial(m_u)
\leq 2$. If $\partial(m_u)=1$ then ${m_u(X)=X-z_0, z_0 \in \C}$,
and therefore  $u=z_0$. Note that the characteristic polynomial is
  $\chi_u(X)=X^2-\textrm{trace}([u])X+\textrm{det}([u]),$  where
trace$([u])=u_1+u_ii+\sigma(u_1+u_ii)=2u_1$ and
 $\textrm{det}([u])=\pm 1$: $$\chi_u(X)=X^2-2u_1X \pm 1.$$
 \para
 \begin{proposition}\label{torp} Let $u=u_1+u_ii+u_jj+u_kk \in \U(\h(\mathfrak o_K))$. Then the following statements hold:\\
\begin{enumerate}
\item $u^2=2u_1u-N(u)$.\\
\item If $N(u)=1$, then $u$ is a torsion unit if, and only if,
 $u_1\in\{-1,0,1\}$ and  the order of $u$ is $1,\,2$, or $4$.\\
\item If $N(u)=-1$, then order of $u$ is  infinite.
\end{enumerate}
\end{proposition}\para
\begin{prf} (1) is obvious.
\para\noindent
(2) Suppose $N(u)=1 $ and $u$ is a torsion unit of order $n$, say. 
If  $X^2-2u_1X+\n(u)=(X-\zeta_1)(X-\zeta_2)$, then $\zeta_i,\
i=1,\,2$, are roots of unity and $\zeta_1\zeta_2=1$.   It follows
that
 $2u_1=\zeta_1+\zeta_2$ is a real number.
 Since $u\in \h(\mathfrak o_K)$ and $\{1,\,\vartheta \}$ is an integral basis of $\mathfrak o_K$, it follows that  $u_1 \in \Z$. From
the equality $2u_1=\zeta_1+\zeta_2$, we have $2|u_1|=
|\zeta_l+\zeta_2|\leq 2$, and therefore $u_1 \in \{-1,\,0,\,1\}$. If
$u_1=0$, then  $u^2=-1$ and therefore
 $o(u)=4$. If $u_1=\pm 1$,  then $\chi_u(X)=X^2 \mp 2X+1=(X \mp 1)^2$, and therefore  $0=\chi_u(u)=(u \mp 1)^2 \in \h(K)$; hence  $u=\pm 1$.
\para\noindent (3)
 If $N(u)=-1$, then $u^2= 2u_1u+1,\ (u^2)_1=2u_1^2+1,\ \n(u^2)=1$. If  $u$  were a torsion unit, then, by  $(2)$ above, $(u^2)_1\in \{-1,\,0.\,1\}$. If $(u^2)_1=0$, then $1/2=-u_1^2\in \mathfrak o_K$, which is not possible.
 If $(u^2)_1=1$, then $u_1=0$, and therefore $u^2=1$ yielding $u=\pm 1$ which is  not the case, because $N(u)=-1$.
 Finally, if $(u^2)_1=-1$, then $u_1^2=-1$ which implies that $\sqrt{-1} \in K$ which is also not the case, because $\h(K)$ is a division ring.
  Hence $u\in \mathcal U(\mathfrak o_K)$  is an element of infinite order.
\end{prf}
\para
 Let $\xi \neq \psi$ be elements of $\{1,\,i,\,j,\,k\}$.
 Suppose \begin{equation}u:=m\sqrt{-d} \xi+p \psi,\  p,\,m\in \mathbb Z,\end{equation}
 is an element in $\h(\mathfrak o_K)$ having  norm $1$.
 Then \begin{equation}\label{Pell}p^2-m^2d=1,\end{equation} i.e., $(p,\,m)$ is a solution of
 the Pell's equation $X^2-dY^2=1$.
 Let $\q:=\mathbb Q[\sqrt{d}]$.  Equation (\ref{Pell}) implies that $\epsilon = p+m\sqrt{d}$
 is a unit in $\mathfrak o_\q$. Conversely, if $\epsilon=p+m\sqrt{d}$ is a unit of norm $1$ in $\mathfrak o_\q$ then, necessarily, $p^2-m^2d=1$, and, therefore,  for any choice of $\xi,\,\psi$ in $\{1,\,i,\,j,\,k\}$, $\xi\neq \psi$, \begin{equation}\label{2pell}
 m\sqrt{-d}\xi+p\psi\end{equation} is a unit in $\h(\mathfrak o_K)$. In particular,  \begin{equation}u_{(\epsilon,\,\psi)}:=p+m\sqrt{-d}\psi,\ \psi\in \{i,\,j,\,k\}, \end{equation} is a unit in $\h(\mathfrak o_K)$.
 \para
 With the notations as above, we have:
 \para
\begin{proposition}\label{gcpell}
\begin{enumerate}
\item If $1 \notin supp(u)$, the support of  $u$, then $u$ is a torsion unit.
\item If  $\epsilon=p + m\sqrt{d}$  is a unit in $\mathfrak o_\q$ then $$
 u^n_{(\epsilon,\,\psi)}=u_{(\epsilon^n,\,\psi)}$$for all $\psi\in \{i,\,j,\,k\} $ and  $n\in \mathbb Z$.
\end{enumerate}
\end{proposition}\para
\begin{prf}

If $1 \notin supp(u)$,  then $u_1=0$; therefore, by the Proposition
\ref{torp},
 $u$ is torsion unit.
 \para
Let  $\mu=A+B\sqrt{d}$ and $\nu=C+D\sqrt{d}$,  be units in $\mathfrak o_\q$. Then
 $u_{(\mu,\psi)}=A  +B \sqrt{-d}\psi$ and
$u_{(\nu,\psi)}=C  +D \sqrt{-d} \psi$ are units in $\h(\mathfrak o_K)$.
We have   $$\mu
  \nu=AC+dBD+(AD+BC)\sqrt{d}.$$ Also  $ u_{(\mu,\psi)} u_{(\nu,\psi)} =(AC+dBD) +(AD+BC)\sqrt{-d}\psi=u_{(\mu \nu,\psi)}$. It follows that  we have $ u^n_{(\epsilon,\,\psi)}=u_{(\epsilon^n,\,\psi)}$ for all $\psi\in \{i,\,j,\,k\} $ and  $n\in \mathbb Z$.
\end{prf}
\para
The units  (\ref{2pell}) constructed above are called $2$-{\it Pell
units}.
\para
\begin{proposition}
Let  $L=\Q[\sqrt{2d}]$, $2d$ square-free,  $\xi, \,\psi,\, \phi$
pairwise distinct elements in $ \{1,\,i,\,j,\,k\}$ and $p,\, m\in
\mathbb Z$. Then the following are equivalent:

\newcounter{3pell}
\begin{list}{$(\roman{3pell})$}{\usecounter{3pell}}

\item
$u:=m\sqrt{-d} \xi + p \psi + (1-p) \phi \in \U(\h(\mathfrak o_K))$.\\
\item
 $\epsilon:=(2p-1)+m\sqrt{2d} \in \U(\mathfrak o_L)$.
 \end{list}
\end{proposition}
\para
\begin{prf}
If $u$ is a unit in $\h(\mathfrak o_K)$ then
$N(u)=-m^2d+p^2+(1-p)^2=1$,  i.e.,  $2p^2-2p-m^2d=0$, and
thus $(2p-1)^2-m^2 2d=1$.  Consequently,
$\epsilon=(2p-1)+m\sqrt{2d}$
 is   invertible in $\mathfrak o_L$. The steps being reversible, the equivalence of (i) and (ii) follows.
\end{prf}
\para
The units constructed above are called $3$-{\it Pell units}. We
shall next determine units of the form
$u=m\sqrt{-d}+(m\sqrt{-d})i+pj+qk$, with $m,\,p,\,q \in \Z$ and
$N(u)=-2m^{2}d+p^{2}+q^{2}=1$. Set $p+q=:r$ and consider the
equation
\begin{equation}\label{4pell}
2p^{2}-2pr-2m^{2}d+r^{2}-1=0.
\end{equation}
\para
\begin{theorem}  If $r=1$, then equation $(\ref{4pell})$
has a solution in $\Z$, and for each such solution, \linebreak
$u=m\sqrt{-d}+(m\sqrt{-d})i+pj+qk$ is a unit in $\h(\mathfrak o_K)$
of norm $1$.
\end{theorem}
\para
\begin{prf} Viewed as a quadratic equation in $p$, (\ref{4pell}) has  real  roots   $$p=\frac{1 \pm
 \sqrt{1+4m^{2}d}}{2}.$$ To obtain a solution in $\Z$, we need the argument under the radical to be a square;
 we thus need to solve the diophantine
 equation
 \begin{equation}\label{5pell}X^{2}-4dY^2=1.\end{equation}

Let $\epsilon=x +y \sqrt{d}$, with $x,\ y\in \mathbb Z$, be a unit in
$\mathfrak o_\q$ having infinite order. Replacing $\epsilon$ by
$\epsilon^2$, if necessary, we can assume that $y$ is even.  We then
have  $x^{2}-y^{2}d=1$, and so $x$ must be  odd.      Taking $m=y/2$
and
  $p=\frac{1 \pm x}{2}$, we obtain a solution of $(\ref{5pell})$ in $\Z$. Clearly,
  for such a solution, the element $u$ lies in $\h(\mathfrak o_K)$ and has norm
  1.
\end{prf}
\para
Using Gauss' result which states that  a positive integer $n$ is a
sum of three  squares if, and only if, $n$ is not of the form
$4^a(8b-1)$, where $a \geq 0$  and  $b \in \mathbb Z$, it is easy to see that, for every integer $m
\equiv 2 \pmod 4$, the integers $m^2 d-1$ and $m^2d+1$ can be expressed as
sums of three  squares.
We can thus construct units $u=m
\sqrt{-d}+pi+qj+rk\in \h(\mathfrak o_K)$ having prescribed norm 1 or
$-1$; we  call such units {\it Gauss units}.
\para
\noindent
{\bf Example.} In   \cite{cjlr}, all units exhibited in $\h(\mathfrak o_{\mathbb
Q[\sqrt{-7}]})$  are of norm $1$. We present some units of norm $-1$
in this ring. The previous theorem guarantees the existence of
integers $p,\,q,\,r$, such that
$$u=6\sqrt{-7}+pi+qj+rk$$ is a unit of norm $-1$. Indeed,
$$(p,\,q,\,r)\in \{(\pm 15,\,\pm5,\,\pm1),\ (\pm13,\,\pm9,\,\pm1),\ (\pm11,\,\pm11,\,\pm3)\},$$
and the triples obtained by permutation of coordinates, are all possible integral solutions.
In \cite{cjlr}, the authors have constructed a   set $S$ of
generators of  the group $SL_1(\h(\mathfrak  o_{\mathbb
Q[\sqrt{-7}]}))$.  If $v_0$ is a unit of  $\h(\mathfrak o_{\mathbb
Q[\sqrt{-7}]})$ having  norm $-1$, then  clearly  $\langle v_0,S
\rangle= \U(\h(\mathfrak o_{\mathbb Q[\sqrt{-7}]}))$.  Thus, for
example, taking ${v_0=6\sqrt{-7}+15i+5j+k}$, we have
\begin{equation}
 \U(\h(\mathfrak o_{\mathbb Q[\sqrt{-7}]}))=\langle v_0,S \rangle.\end{equation}
\para

The set $\{1,\frac{1+\sqrt{-7}}{2}\}$ is an integral basis of  $R=\Z[\frac{1+\sqrt{-7}}{2}]$.
 Consider units of the form $$\frac{m+\sqrt{-d}}{2}\pm(\frac{m-\sqrt{-d}}{2})i+pj$$
These are  neither  Pell nor   Gauss units. Those of norm $\pm 1$, are solutions of the equation
\begin{equation}\label{eqgr}m^2+2p^2=\pm 2+d \end{equation} in $\Z$.
The main result of \cite{cjlr} states that if  $d=7$ then the units of norm $1$ of  the above type, together  with
the trivial units $i$ and $j$, generate the group $SL_1(\h(R))$.

For $d \equiv 7 \pmod 8$ there are no units of norm $-1$ of the above type, since, in this case,
the equation $m^2+2p^2=-2+d$  has no solution in $\Z$, as can be easily seen working module $8$.

\typeout{In case the unit group is generated by units satisfying diophantine quadratic equations.
So our construction might indicate an algebraic way of obtaining generators for the
unit group $\mathcal U_1(\mathfrak o_K[Q_8])$ in the general case. }
In case $d\neq 7$, we give some more examples  of negative norm units of the form \linebreak
$\frac{m+\sqrt{-d}}{2}\pm(\frac{m-\sqrt{-d}}{2})i+pj$.

If $d=15$ then the equation (\ref{eqgr}) becomes $m^2+2p^2=17$; the pairs 
$(m,p)\in \{(3,2),(3,-2),\linebreak (-3,2),(-3,-2)\}$ are its integral solutions. For $m=3$ either $p=2$ or
$p=-2$ and so there are $8$ units. Each coefficient of $u$ is distinct, hence
for each solution $(m,p)$ there are
$3!$ units with the same support, thus there are $36$ different units for a given fixed support. By Proposition \ref{torp} all these units have
 infinite order if $u_1 \notin \{-1,0,1\}$. If $1 \in supp(u)$ then either $\{i,j\}
\subset supp(u)$ or $\{i,k\} \subset supp(u)$, or $\{j,k\} \subset supp(u)$. Therefore there are
 $108$ of these units and, for example,    $$\frac{3+\sqrt{-15}}{2}+(\frac{3-\sqrt{-15}}{2})j-2k$$ is one of them.

If $1 \notin
supp(u)$ then  $u$ is a torsion unit, so there are $36$ torsion units of this type.
 One of them is the unit $$(\frac{-3-\sqrt{-15}}{2})i+(\frac{-3+\sqrt{-15}}{2})j+2k,$$ of order $4$.

For $d=31$ we obtain  $m^2+2p^2=33$ whose  solutions in $\Z$ are:
$(m,p)\in \{(1,4),(1,-4),(-1,4),\linebreak (-1,-4)\}$.

As another example of a unit of norm $-1$ in a quaternion algebra,
we may mention that, in  $\h(\mathfrak o_{\mathbb Q[\sqrt{-23}]})$,
$u=5\sqrt{-23}+23i+6j+3k$ is a unit of norm $-1$.

\para
We next exhibit some Gauss units of norm $1$. For $\h(\mathfrak
o_{\mathbb Q[\sqrt{-15}]})$, there exist $p,\,q,\,r$, such that
$u=10\sqrt{-15}+pi+qj+rk$ is a unit of norm $1$. In fact, 
$(36,\,14,\,3),\ (36,\,13,\,6),$ $\ (32,\,21,\,6),\ (30,\,24,\,5)$ are
some  of the possible choices  for $(p,\,q,\,r)$. For $\h(\mathfrak
o_{\mathbb Q[\sqrt{-23}]})$, $u=2\sqrt{-23}+8i+5j+2k$ is a unit of
norm $1$. It is interesting to note that 
${u=3588\sqrt{-23}+12168i+12167j}$ is a Gauss unit, although $4$ divides
$3588$.
\para

\para

We conclude with the following result: \para
\begin{theorem}
Let $K=\Q [\sqrt{-d}]$, $0<d\equiv 7 \pmod 8$ and $\mathfrak o_K$ the ring of integers of $K$. If $\epsilon=p+m\sqrt{d}$  is a
unit in $\Z [\sqrt{d}]$, and
$x:=u_{(\epsilon,\,\psi)},\ y:=u_{(\epsilon,\,\psi')}$ are two   $2$-Pell
units in $\mathcal U(\h(\mathfrak o_K))$, where $ \psi$ and
$\psi'\in \{i,\,j,\,k\}$ and $\psi\neq \psi'$, then there exists a
natural number $m$ such that $\langle x^m,\,y^m \rangle$ is a free
group of rank $2$.
\end{theorem}
\para
\begin{prf}
By Proposition \ref{hypq8}, $\U(\h(\mathfrak o_K))$ is a hyperbolic
group. In view of \cite[Proposition $III.\Gamma.3.20$]{mbah}, there
exists a natural $m$, such that,
 $\langle x^m,\,y^m \rangle$ is a free group of rank at most $2$. However,  Proposition \ref{gcpell}  item
$(2)$ ensures that  $\langle x\rangle\cap \langle y\rangle=\{1\}$.
Therefore,
  $\langle x^m,\,y^m \rangle$ has rank at least 2, and hence  $\langle x^m,y^m \rangle$
 is a free group of rank $2$.
\end{prf}

\par

\vspace{3cm}

\begin{flushleft}Instituto de Matem\'atica e Estat\'\i stica,\\
 Universidade de S\~ao Paulo (IME-USP),\\ Caixa Postal 66281, S\~ao Paulo,\\ CEP
 05315-970 - Brasil \\ email  -- ostanley@ime.usp.br -- \end{flushleft} 
 \begin{flushleft} Centre for Advanced Study in Mathematics, \\Panjab University, \\Chandigarh 160014 - India \\ email -- ibspassi@yahoo.co.in --
 \end{flushleft}
 \begin{flushleft}Escola de Artes, Ci\^encias e Humanidades,\\
 Universidade de S\~ao Paulo (EACH-USP),\\ Rua Arlindo B\'ettio, 1000, Ermelindo Matarazzo, S\~ao Paulo, \\CEP
 03828-000 - Brasil \\email  -- acsouzafilho@usp.br -- \end{flushleft}

\end{document}